\documentclass[12pt]{amsart}

\usepackage{amssymb,amscd,psfig,latexsym,graphicx,eucal}

\usepackage{fullpage}

\newtheorem{theorem}{Theorem}[section]
\newtheorem{lemma}[theorem]{Lemma}

\theoremstyle{definition}

\numberwithin{equation}{section}

\newcommand{\diam}{\operatorname{diam}}
\newcommand{\area}{\operatorname{area}}

\newcommand{\ds}{\displaystyle}

\newcommand{\Hdim}{\operatorname{dim}_{\operatorname{H}}}

\newcommand{\ve}{\varepsilon}

\newcommand{\sm}{\smallsetminus}
\newcommand{\bd}{\partial}
\newcommand{\Chat}{\widehat{\mathbb C}}

\newcommand{\ov}{\overline}

\newcommand{\bl}{\boldsymbol{d}}

\newcommand{\CC}{{\mathbb C}}

\newcommand{\RR}{{\mathbb R}}

\newcommand{\DD}{{\mathbb D}}

\newcommand{\vs}{\vspace{2mm}}

\newcommand{\Sen}{{\mathbb  S}^1}

\renewcommand{\marginpar}[1]{}
\catcode`\@=12

\def\Empty{}
\newcommand\oplabel[1]{
  \def\OpArg{#1} \ifx \OpArg\Empty {} \else
  	\label{#1}
  \fi}
		
\long\def\realfig#1#2#3#4{
\begin{figure}[tp]
\centerline{\psfig{figure=#2,width=#4}}
\caption[#1]{#3}
\oplabel{#1}
\end{figure}}

\newcommand{\comm}[1]{}

\newcommand{\cal}{\mathcal}

\newcommand{\lemref}[1]{Lemma~\ref{#1}}

\newcommand{\figref}[1]{Fig.~\ref{#1}}

\begin{document}

\title
{David maps and Hausdorff Dimension}

\author[S. Zakeri]
{Saeed Zakeri}

\address{S. Zakeri, Institute for Mathematical Sciences, Stony Brook
University, Stony Brook, NY 11794-3651}
\email{zakeri@math.sunysb.edu}



\date{November 19, 2002}

\begin{abstract}
David maps are generalizations of classical planar
quasiconformal maps for which the
dilatation is allowed to tend to infinity in a
controlled fashion. In this note we examine how these
maps distort Hausdorff dimension. We show \vs
\begin{enumerate}
\item[$\bullet$]
Given $\alpha$ and $\beta$ in $[0,2]$, there exists a David
map  $\varphi:\CC \to \CC$ and a compact set $\Lambda$ such that
$\Hdim \Lambda =\alpha$ and $\Hdim \varphi(\Lambda)=\beta$. \vs
\item[$\bullet$]
There exists a David map $\varphi:\CC \to \CC$ such that
the Jordan curve $\Gamma=\varphi (\Sen)$ satisfies $\Hdim
\Gamma=2$.\vs
\end{enumerate}
One should contrast the first statement with the fact that
quasiconformal maps preserve sets of Hausdorff dimension $0$ and
$2$. The second statement provides an example of a Jordan curve
with Hausdorff dimension $2$ which is (quasi)conformally
removable.
\end{abstract}

\maketitle
\thispagestyle{empty} \def\IMSmarkvadjust{0 pt}
\def\IMSmarkhadjust{0 pt}
\def\IMSmarkhpadding{0 pt}
\def\IMSpubltext{Published in modified form:}
\def\SBIMSMark#1#2#3{
 \font\SBF=cmss10 at 10 true pt
 \font\SBI=cmssi10 at 10 true pt
 \setbox0=\hbox{\SBF \hbox to \IMSmarkhpadding{\relax}
                Stony Brook IMS Preprint \##1}
 \setbox2=\hbox to \wd0{\hfil \SBI #2}
 \setbox4=\hbox to \wd0{\hfil \SBI #3}
 \setbox6=\hbox to \wd0{\hss
             \vbox{\hsize=\wd0 \parskip=0pt \baselineskip=10 true pt
                   \copy0 \break%
                   \copy2 \break%
                   \copy4 \break}}
 \dimen0=\ht6   \advance\dimen0 by \vsize \advance\dimen0 by 8 true pt
                \advance\dimen0 by -\pagetotal
	        \advance\dimen0 by \IMSmarkvadjust
 \dimen2=\hsize \advance\dimen2 by .25 true in
	        \advance\dimen2 by \IMSmarkhadjust

%
%
  \openin2=publishd.tex
  \ifeof2\setbox0=\hbox to 0pt{}
  \else 
     \setbox0=\hbox to 3.1 true in{
                \vbox to \ht6{\hsize=3 true in \parskip=0pt  \noindent  
                {\SBI \IMSpubltext}\hfil\break
                \input publishd.tex 
                \vfill}}
  \fi
  \closein2
  \ht0=0pt \dp0=0pt
 \ht6=0pt \dp6=0pt
 \setbox8=\vbox to \dimen0{\vfill \hbox to \dimen2{\copy0 \hss \copy6}}
 \ht8=0pt \dp8=0pt \wd8=0pt
 \copy8
 \message{*** Stony Brook IMS Preprint #1, #2. #3 ***}
}

\def\IMSmarkvadjust{-30pt}
\SBIMSMark{2002/05}{November 2002}{}

\section{Introduction}
\label{sec:intro}

An orientation-preserving homeomorphism $\varphi: U \to V$ between
planar domains is called {\it quasiconformal} if it belongs to the
Sobolev class $W^{1,1}_{\text{loc}}(U)$ (i.e., has locally
integrable distributional partial derivatives in $U$) and its
complex dilatation $\mu_{\varphi}:=\ov{\bd}\varphi / \bd \varphi$
satisfies
$$\| \mu_{\varphi} \|_{\infty} <1.$$
In terms of the real dilatation defined by
$$K_{\varphi}:=\frac{1+|\mu_{\varphi}|}{1-|\mu_{\varphi}|}=
\frac{|\bd \varphi|+|\ov{\bd}\varphi|}{|\bd \varphi|-|\ov{\bd}\varphi|},$$
the latter condition can be expressed as
$$\| K_{\varphi} \|_{\infty} < +\infty.$$
The quantity $\| K_{\varphi} \|_{\infty}$ is called the {\it
maximal dilatation} of $\varphi$. We say that $\varphi$ is
$K$-quasiconformal if its maximal dilatation does not exceed $K$.

For later comparison with the properties of David maps
defined below, we recall some basic properties of quasiconformal
maps (see \cite{Ahlfors} or \cite{LV}):
\begin{enumerate}
\item[$\bullet$]
If $\varphi$ is $K$-quasiconformal for some $K \geq 1$, so
is the inverse map $\varphi^{-1}$. \vs
\item[$\bullet$]
A $K$-quasiconformal map $\varphi: U \to V$ is locally H\"{o}lder
continuous of exponent $1/K$. In other words, for every compact set
$E \subset U$ and every $z,w \in E$,
$$|\varphi(z)-\varphi(w)| \leq C\, |z-w|^{\frac{1}{K}}$$
where $C>0$ only depends on $E$ and $K$. \vs
\item[$\bullet$]
A quasiconformal map $\varphi: U \to V$ is absolutely continuous;
in fact, the Jacobian $J_{\varphi}=|\bd \varphi|^2 - |\ov{\bd}
\varphi|^2$ is locally integrable in $U$ and
\begin{equation}\label{eqn:JJ}
\area \varphi(E)=\int_E J_{\varphi} \, dx\, dy,
\end{equation}
for every measurable $E \subset U$. \vs
\item[$\bullet$]
More precisely, the Jacobian $J_{\varphi}$ of a quasiconformal map
$\varphi: U \to V$ is in $L^p_{\text{loc}}(U)$ for some $p>1$. If
we define
\begin{equation}
\label{eqn:pk}
p(K):= \sup \{ p: J_{\varphi} \in L^p_{\text{loc}}(U) \ \text{for
every} \ K\text{-quasiconformal map} \ \varphi \ \text{in} \ U \},
\end{equation}
then $p(K)$ is independent of the domain $U$ and
\begin{equation}
\label{eqn:vpk}
p(K)=\frac{K}{K-1}.
\end{equation}
This was conjectured by Gehring and V\"ais\"al\"a in 1971 \cite{GV}
and was proved by Astala in 1994 \cite{Astala}. \vs
\item[$\bullet$]
Let $\{ \varphi_n \}$ be a sequence of $K$-quasiconformal maps in a planar
domain $U$ which fix two given points of $U$. Then $\{ \varphi_n
\}$ has a subsequence which converges locally uniformly to a
$K$-quasiconformal map in $U$.
\end{enumerate}

The measurable Riemann mapping theorem of Morrey-Ahlfors-Bers
\cite{AB} asserts that any measurable function $\mu$ in a domain
$U$ which satisfies $\| \mu \|_{\infty} <1$ is the complex
dilatation of some quasiconformal map $\varphi$ in $U$, which means
$\varphi$ satisfies the {\it Beltrami equation} $\ov{\bd} \varphi =
\mu \cdot \bd \varphi$ almost everywhere in $U$. Recent progress in
conformal geometry and holomorphic dynamics has made it abundantly
clear that one must also study this equation in the case $\| \mu
\|_{\infty} =1$. With some restrictions on the asymptotic growth of
$|\mu|$, the solvability of the Beltrami equation can still be
guaranteed. One such condition is given by David in \cite{David}.
Let $\sigma$ denote the spherical area in $\Chat$ and $\mu$ be a
measurable function in $U$ which satisfies
\begin{equation}
\label{eqn:David1}
\sigma \{ z \in U : |\mu(z)| > 1-\ve \} \leq C \exp \left(
-\frac{\alpha}{\ve} \right) \qquad \text{for all}\ \ve < \ve_0
\end{equation}
for some positive constants $C, \alpha, \ve_0$. Then David showed that
the Beltrami equation $\ov{\bd} \varphi = \mu \cdot \bd \varphi$
has a homeomorphic solution $\varphi \in W^{1,1}_{\text{loc}}(U)$
which is unique up to postcomposition with a
conformal map. Motivated by this result, we call a homeomorphism
$\varphi: U \to V$ a {\it David map} if $\varphi \in
W^{1,1}_{\text{loc}}(U)$ and the complex dilatation $\mu_{\varphi}$
satisfies a condition of the form \eqref{eqn:David1}.
Equivalently, $\varphi$ is a David map if
there are positive constants $C, \alpha, K_0$ such that its real
dilatation satisfies
\begin{equation}
\label{eqn:David2}
\sigma \{ z \in U : K_{\varphi}(z) > K \} \leq C e^{-\alpha K} \qquad
\text{for all}\ K>K_0.
\end{equation}
To emphasize the values of these constants, sometimes we say that
$\varphi$ is a $(C, \alpha, K_0)$-David map. Note that when
$U$ is a bounded domain in $\CC$, the spherical metric in
\eqref{eqn:David1} or \eqref{eqn:David2} can be replaced with the
Euclidean area.

David maps enjoy some of the useful properties of
quasiconformal maps, but the two classes differ in many respects.
As indications of their similarity, let us mention the following
two facts:
\begin{enumerate}
\item[$\bullet$]
Every David map is absolutely continuous; the Jacobian
formula \eqref{eqn:JJ} still holds.\vs
\item[$\bullet$]
{\it Tukia's Theorem} \cite{Tukia}. ``Let $C,\alpha,K_0$ be positive and
suppose $\{ \varphi_n \}$ is a sequence of
$(C,\alpha,K_0)$-David maps in a domain $U$ which fix two given
points of $U$. Then $\{ \varphi_n \}$ has a subsequence which
converges locally uniformly to a David map in $U$.'' It
is rather easy to show that some subsequence of $\{ \varphi_n \}$
converges locally uniformly to a homeomorphism, but that this
homeomorphism must be David is quite non-trivial. We remark that
the parameters of the limit map may a priori be different from
$C,\alpha,K_0$.
\end{enumerate}

Here are further properties of David maps which indicate
their difference with quasiconformal maps:
\begin{enumerate}
\item[$\bullet$]
The inverse of a David map may not be David.\vs
\item[$\bullet$]
A David map may not be locally H\"{o}lder. \vs
\item[$\bullet$]
The Jacobian of a David map may not be in $L^p_{\text{loc}}(U)$ for
any $p>1$.
\end{enumerate}
As an example, the homeomorphism $\varphi: \DD(0,e^{-1}) \to \DD$
defined by
$$\varphi(re^{i\theta}):=-\frac{1}{\log r}\, e^{i\theta}$$
is a David map but $\varphi^{-1}$ is not. Moreover, $\varphi$ is not
H\"{o}lder in any neighborhood of $0$, and $J_{\varphi} \notin
L^p_{\text{loc}}$ for $p>1$. \vs

The main goal of this note is to show how David maps
differ from quasiconformal maps in the way they change Hausdorff
dimension of sets. Recall that the {\it Hausdorff $s$-measure} of
$E \subset \CC$ is defined by
$$H^s(E):= \lim_{\ve \to 0} \inf_{\cal U} \sum_i (\diam U_i)^s,$$
where the infimum is taken over all countable covers ${\cal U}=\{
U_i \}$ of $E$ by sets of Euclidean diameter at most
$\ve$. The {\it Hausdorff dimension} of $E$ is defined by
$$\Hdim E := \inf \{ s : H^s(E)=0 \}.$$
Quasiconformal maps can change Hausdorff dimension of sets only by
a bounded factor depending on their maximal dilatation. This was
first proved by Gehring and V\"ais\"al\"a \cite{GV} who showed that
if $\varphi:U \to V$ is $K$-quasiconformal, $E \subset U$, $\Hdim
E=\alpha$ and $\Hdim \varphi(E)=\beta$, then
$$\frac{2(p(K)-1)\alpha}{2p(K)-\alpha} \leq \beta \leq \frac{2 p(K)
\alpha}{2(p(K)-1)+\alpha}.$$
Here $p(K)>1$ is the constant defined
in \eqref{eqn:pk}. By Astala's result \eqref{eqn:vpk}, one obtains
$$\frac{2\alpha}{2K-(K-1)\alpha} \leq \beta \leq
\frac{2K\alpha}{2+(K-1)\alpha}$$ which can be put in the symmetric
form
\begin{equation}
\label{eqn:sym}
\frac{1}{K} \left( \frac{1}{\alpha} - \frac{1}{2} \right) \leq
\frac{1}{\beta} - \frac{1}{2} \leq K \left( \frac{1}{\alpha} -
\frac{1}{2} \right).
\end{equation}
It follows in particular that quasiconformal maps preserve sets of
Hausdorff dimension $0$ and $2$.

By contrast, we prove \vs \\
{\bf Theorem A.} {\it Given any two numbers $\alpha$
and $\beta$ in $[0,2]$, there exists a David map $\varphi:\CC
\to \CC$ and a compact set $\Lambda \subset \CC$ such that
$\Hdim \Lambda=\alpha$ and $\Hdim \varphi(\Lambda)=\beta$.} \vs \\
The proof shows that the parameters of $\varphi$ can be taken 
independent of $\alpha$ and $\beta$.

In the special case of a {\it $K$-quasicircle}, i.e., the image
$\Gamma$ of the round circle under a $K$-quasiconformal map, the estimate
\eqref{eqn:sym} gives
$$1 \leq \Hdim \Gamma \leq \frac{2K}{K+1}$$
(the lower bound comes from topological considerations). It is
well-known that $\Hdim \Gamma$ can in fact take all values in $[1,2)$.
We show that the upper bound $2$ is attained by a David image of
the round circle. Let us call a Jordan curve $\Gamma \subset \CC$
a {\it David circle} if there exists a David map $\varphi: \CC \to \CC$
such that $\Gamma=\varphi(\Sen)$, where
$\Sen$ is the unit circle $\{ z \in \CC: |z|=1 \}$. \vs \\
{\bf Theorem B.} {\it There exist David circles of Hausdorff
dimension $2$.}\vs \\
One corollary of this result is that there are Jordan curves of
Hausdorff dimension $2$ that are (quasi)conformally removable (see
\S \ref{sec:proofB}).

Both results are bad (or exciting?) news for applications in
holomorphic dynamics, where one often wants to estimate the Hausdorff 
dimension of invariant sets by computing the dimension in a conjugate 
dynamical system. The dichotomy of having dimension $<2$ or $=2$ for such
invariant sets, which is respected by quasiconformal conjugacies,
is no longer preserved by David conjugacies.
For example, by performing quasiconformal surgery on a Blaschke
product, Petersen proved that the Julia set of the quadratic polynomial 
$Q_{\theta}: z \mapsto e^{2 \pi i \theta}z+z^2$ is locally-connected and 
has measure zero whenever $\theta$ is an irrational of bounded type 
\cite{Petersen}. In this case, the boundary of the Siegel disk of 
$Q_{\theta}$ is a quasicircle whose Hausdorff dimension is strictly 
between $1$ and $2$ (compare \cite{GJ}). On the other hand, by performing a 
{\it trans-quasiconformal} surgery and using David's theorem, 
Petersen and the author extended the above result to almost
every $\theta$ \cite{PZ}. It follows that there exists a full-measure set of
rotation numbers $\theta$ for which the boundary of the Siegel disk of 
$Q_{\theta}$ is a David circle but not a quasicircle. Thus, Theorem B opens  
the possibility that this boundary alone might have dimension $2$,
which would be a rather curious phenomenon.   

\section{Preliminary constructions}
\label{sec:bg}

For two positive numbers $a$ and $b$, we write
$$a \preccurlyeq b$$
if there is a universal constant $C>0$ such that $a \leq C b$. We
write
$$a \asymp b$$
if $a \preccurlyeq b$ and $b \preccurlyeq a$, i.e., if there is a
universal constant $C>0$ such that $C^{-1}b \leq a \leq C b$. In this
case, we say that $a$ and $b$ are {\it comparable}.

\subsection*{A family of Cantor sets}
Given a strictly decreasing sequence $\bl = \{ d_n \}_{n \geq 0}$
of positive numbers with $d_0=1$, we construct a Cantor set
$\Lambda(\bl)$ as the intersection of a nested sequence $\{
\Lambda_n \}_{n \geq 0}$ of compact sets in the unit square
$\Lambda_0:=[-\frac{1}{2},\frac{1}{2}] \times
[-\frac{1}{2},\frac{1}{2}]$ defined inductively as follows. Set
$a_1:=2^{-2}(d_0-d_1)$ and define $\Lambda_1$ as the disjoint union
of the four closed squares of side-length $2^{-1}d_1$ in
$\Lambda_0$ which have distance $a_1$ to the boundary of
$\Lambda_0$ (see \figref{cantor}). Suppose $\Lambda_{n-1}$ is
constructed for some $n \geq 2$ so that it is the disjoint union
of $4^{n-1}$ closed squares of side-length $2^{-(n-1)}d_{n-1}$.
Define
\begin{equation}
\label{eqn:a_n}
a_n:=2^{-(n+1)}(d_{n-1}-d_n).
\end{equation}
For any square $S$ in $\Lambda_{n-1}$, consider the disjoint union
of the four closed squares in $S$ of side-length $2^{-n}d_n$ which
have distance $a_n$ to the boundary of $S$. The union of all these
squares for all such $S$ will then be called $\Lambda_n$. Clearly
$\Lambda_n$ is the disjoint union of $4^n$ closed squares of
side-length $2^{-n}d_{n}$, and the inductive definition is
complete.

\realfig{cantor}{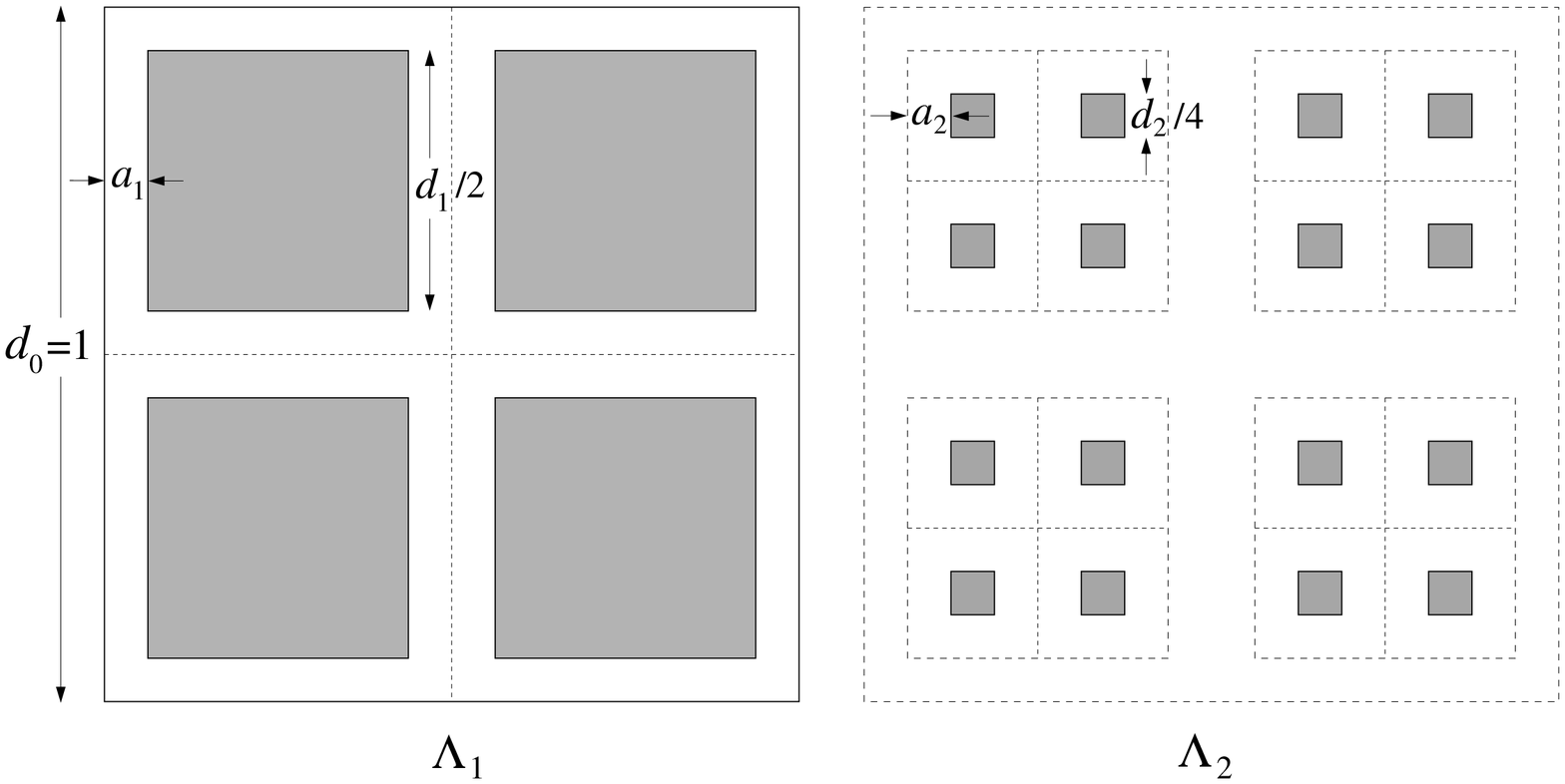}{{\sl First two steps in the construction
of $\Lambda(\bl)$.}}{12cm}

The Cantor set $\Lambda(\bl)$ is defined as $\bigcap_{n\geq 0}
\Lambda_n$. We have
$$\area \Lambda(\bl) = \lim_{n \to \infty} \area \Lambda_n = \lim_{n
\to \infty} d_n^2.$$

\begin{lemma}
\label{Hdim}
The Hausdorff dimension of the Cantor set $\Lambda=\Lambda(\bl)$
satisfies
\begin{equation}
\label{eqn:dimest}
2-\limsup_{n \to \infty} \frac{-2\log d_{n+1}}{-\log d_n + n \log
2} \leq \Hdim \Lambda \leq 2-\liminf_{n \to \infty} \frac{-2\log
d_n}{-\log d_n + n \log 2}.
\end{equation}
\end{lemma}

\begin{proof}
For each $n\geq 0$, there are $4^n$ squares of diameter $2^{\frac{1}{2}-n}
d_n$ covering $\Lambda$. Hence the Hausdorff $s$-measure of
$\Lambda$ is bounded above by
$$\liminf_{n \to \infty}\, 4^n (2^{\frac{1}{2}-n}d_n)^s=2^{\frac{s}{2}}\,
\liminf_{n \to \infty}\, 2^{n(2-s)}d_n^s,$$
which is zero if $s>2-\liminf_{n \to \infty} (-2\log d_n)/(-\log
d_n + n \log 2)$. This proves the upper bound in
(\ref{eqn:dimest}).

The lower bound follows from a standard mass distribution argument:
Construct a probability measure $\mu$ on $\Lambda$ which
gives equal mass $4^{-n}$ to each square in $\Lambda_n$, so that
$$
\mu(S) = \frac{\area(S)}{d_n^2} \quad \text{if}\ S \ \text{is a
square in} \ \Lambda_n.
$$
Let $x \in \Lambda$ and $\ve >0$, and choose $n$ so that $2^{-n}d_n
< \ve \leq 2^{-(n-1)}d_{n-1}$. The disk $\DD(x,\ve)$ intersects at
most $\pi \ve^2 / (4^{-n} d_n^2)$ squares in $\Lambda_n$ each
having $\mu$-mass of $4^{-n}$. It follows that
$$
\mu(\DD(x,\ve)) \preccurlyeq \frac{\ve^2}{d_n^2} =
\ve^s \, \frac{\ve^{2-s}}{d_n^2} \preccurlyeq \ve^s
\, \frac{2^{-n(2-s)}d_{n-1}^{2-s}}{d_n^2}.
$$
If $s<2-\limsup_{n \to \infty} (-2\log d_{n+1})/(-\log d_n + n \log
2)$, the term $2^{-n(2-s)}d_{n-1}^{2-s}/d_n^2$ will tend to zero as $n
\to \infty$, so that
$$
\mu(\DD(x,\ve)) \preccurlyeq \ve^s.
$$
It follows from Frostman's Lemma (see for example \cite{Mattila})
that $\Hdim \Lambda \geq s$. This gives the lower bound in
(\ref{eqn:dimest}).
\end{proof}

\subsection*{Standard homeomorphisms between Cantor sets}
We construct standard homeomorphisms with controlled dilatation
between Cantor sets of the form $\Lambda(\bl)$ defined above. The
construction will depend on the following lemma:

\begin{lemma}
\label{annulus}
Fix $0<a \leq b< \frac{1}{2}$. Let $A_a$ be the closed annulus
bounded by the squares
$$
\left \{ (x,y) \in \RR^2 : \max \{ |x|,|y| \} =\frac{1}{2} \right \}
\quad \text{and} \quad \left \{ (x,y) \in \RR^2 : \max \{ |x|,|y|
\} = \frac{1}{2}-a \right \},
$$
and similarly define $A_b$. Let $\varphi: \bd A_a \to \bd A_b$ be a
homeomorphism which is the identity on the outer boundary component
and acts affinely on the inner boundary component, mapping
$\frac{1}{2}-a$ to $\frac{1}{2}-b$. Then $\varphi$ can be extended
to a $K$-quasiconformal homeomorphism $A_a \to A_b$, with
\begin{equation}
\label{eqn:KK}
K \asymp \frac{b\, (1-2a)}{a\, (1-2b)}.
\end{equation}
\end{lemma}

\begin{proof}
Let us first make a simple observation: If $z$ and $w$ are points
in the upper half-plane and $L: \RR^2 \to \RR^2$ is the affine map
such that $L(0)=0$, $L(1)=1$ and $L(z)=w$ (see \figref{ann}), then
the real dilatation of $L$ is given by
\begin{equation}
\label{eqn:KL}
K_L=\frac{|z-\ov{w}|+|z-w|}{|z-\ov{w}|-|z-w|}.
\end{equation}
To prove the lemma, take the triangulations of $A_a$ and $A_b$
shown in \figref{ann} and extend $\varphi$ affinely to each
triangle. After appropriate rescaling, it follows from
\eqref{eqn:KL} that on a triangle of type I in the figure, the
dilatation of $\varphi$ is comparable to $b/a$, while on a triangle
of type II, the dilatation of $\varphi$ is comparable to
$b(1-2a)/(a(1-2b))$. Since $b(1-2a)/(a(1-2b)) \geq b/a$, we obtain
\eqref{eqn:KK}.
\end{proof}

\realfig{ann}{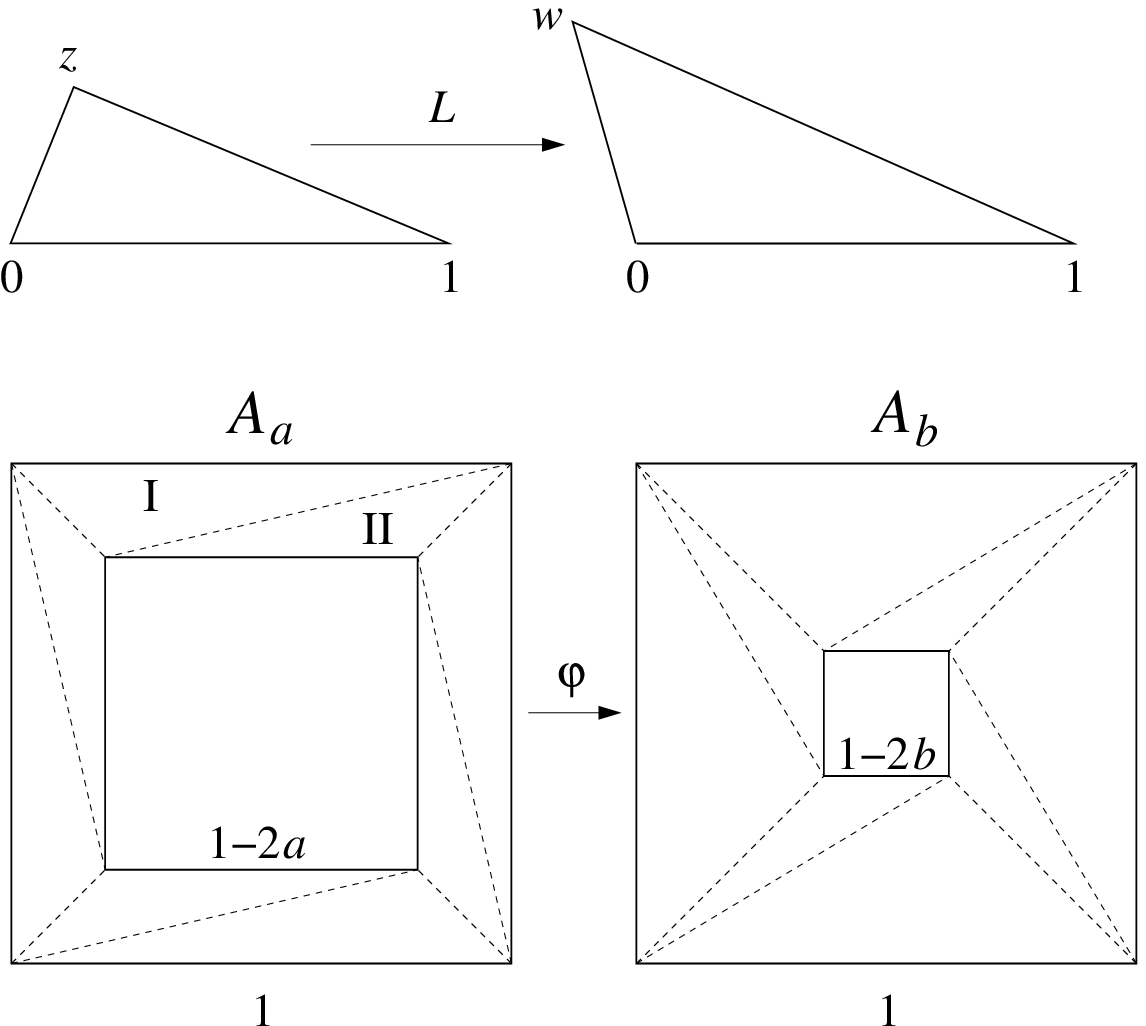}{}{9cm}

Now take a decreasing sequence $\bl= \{ d_n \}$ of positive numbers
with $d_0=1$, let $\{ a_n \}$ be defined as in (\ref{eqn:a_n}), and
consider the Cantor set $\Lambda(\bl) =\bigcap \Lambda_n$. Take
another such sequence $\bl' =\{ d'_n \}$ and let $a'_n, \Lambda'_n,
\Lambda(\bl')$ denote the corresponding data. We construct a
homeomorphism $\varphi: \CC \to \CC$ which maps the Cantor set
$\Lambda=\Lambda(\bl)$ to $\Lambda'=\Lambda(\bl')$. This $\varphi$
is the uniform limit of a sequence of quasiconformal
maps $\varphi_n: \CC \to \CC$ with
$\varphi_n(\Lambda_n)=\Lambda'_n$, defined inductively as follows.
Let $\varphi_0$ be the identity map on $\CC$. Suppose
$\varphi_{n-1}$ is constructed for some $n \geq 1$ and that it maps
each square in $\Lambda_{n-1}$ affinely to the corresponding square
in $\Lambda'_{n-1}$. Define $\varphi_n=\varphi_{n-1}$ on $\CC \sm
\Lambda_{n-1}$ and let $\varphi_n$ map each square in $\Lambda_n$
affinely to the corresponding square in $\Lambda'_n$. The remaining
set $\Lambda_{n-1} \sm \Lambda_n$ is the union of $4^n$ annuli on
the boundary of which $\varphi_n$ can be defined affinely. By
rescaling each annulus in $\Lambda_{n-1} \sm \Lambda_n$ and the
corresponding annulus in $\Lambda'_{n-1} \sm \Lambda'_n$, we are in
the situation of \lemref{annulus}, so we can extend $\varphi_n$ in
a piecewise affine fashion to each such annulus. This defines
$\varphi_n$ everywhere, and the inductive definition is complete.

To estimate the maximal dilatation of $\varphi_n$, note that by the
above construction $\varphi_n$ is conformal in $\Lambda_n$ and has
the same dilatation as $\varphi_{n-1}$ on $\CC \sm \Lambda_{n-1}$.
On each of the $4^n$ annuli in $\Lambda_{n-1} \sm \Lambda_n$, the
dilatation of $\varphi_n$ can be estimated using \eqref{eqn:KK} in
\lemref{annulus}. In fact, rescaling each such annulus by a factor
$2^n/d_{n-1}$ and the corresponding annulus in
$\Lambda'_{n-1} \sm \Lambda'_n$ by a factor $2^n/d'_{n-1}$,
it follows from \eqref{eqn:KK} that the dilatation of $\varphi_n$
on each such annulus is comparable to
\begin{align*}
& \max \left\{
\frac{\frac{a_n'}{2^{-n}d'_{n-1}}(1-2\frac{a_n}{2^{-n}d_{n-1}})}
{\frac{a_n}{2^{-n}d_{n-1}}(1-2\frac{a'_n}{2^{-n}d'_{n-1}})},
\frac{\frac{a_n}{2^{-n}d_{n-1}}(1-2\frac{a'_n}{2^{-n}d'_{n-1}})}
{\frac{a_n'}{2^{-n}d'_{n-1}}(1-2\frac{a_n}{2^{-n}d_{n-1}})}
\right \} \vs \vs \\
= & \max \left\{
\frac{a_n'(d_{n-1}-2^{n+1}a_n)}{a_n(d'_{n-1}-2^{n+1}a'_n)},
\frac{a_n(d'_{n-1}-2^{n+1}a'_n)}{a_n'(d_{n-1}-2^{n+1}a_n)}
\right \} \vs \vs \\
= & \max \left\{
\frac{a'_n d_n}{a_n d'_n},
\frac{a_n d'_n}{a'_n d_n}
\right \} .
\end{align*}
To sum up, the construction gives a sequence $\{ \varphi_n \}$ with the
following properties:
\begin{enumerate}
\item[(i)]
$\varphi_n=\varphi_{n-1}$ on $\CC \sm \Lambda_{n-1}$. \vs
\item[(ii)]
$\varphi_n$ maps each square in $\Lambda_n$ affinely to the
corresponding square in $\Lambda'_n$. \vs
\item[(iii)]
$\varphi_n$ is $K_n$-quasiconformal, where
\begin{equation}
\label{eqn:dil}
K_n \asymp \max \left\{ K_{n-1}, \frac{a_n' d_{n}}{a_n d_{n}'},
\frac{a_n d_{n}'}{a_n' d_{n}} \right \}
\end{equation}
and $K_0=1$.
\end{enumerate}
Evidently, $\varphi:=\lim_{n \to \infty} \varphi_n$ is a
homeomorphism which agrees with $\varphi_n$ on $\CC \sm \Lambda_n$
for every $n$ and
satisfies $\varphi (\Lambda)=\Lambda'$. We call this $\varphi$
the {\it standard homeomorphism from $\Lambda$ to $\Lambda'$}. Observe 
that by the construction, the inverse map $\varphi^{-1}$ is the
standard homeomorphism from $\Lambda'$ to $\Lambda$.

\section{Proof of Theorem A}
\label{sec:proofA}

We are now ready to prove Theorem A cited in \S \ref{sec:intro}. \vs \\
{\it Proof of Theorem A.} \ If $0 < \alpha, \beta < 2$, it is
well-known that there is a $K$-quasiconformal map $\varphi:\CC \to
\CC$ mapping a set of dimension $\alpha$ to a set of dimension
$\beta$ (see for example \cite{GV}). Moreover, by \eqref{eqn:sym},
the minimum $K$ this would require is
$$
\max \left \{ \frac{\frac{1}{\beta}-\frac{1}{2}}
{\frac{1}{\alpha}-\frac{1}{2}}\ , \
\frac{\frac{1}{\alpha}-\frac{1}{2}}
{\frac{1}{\beta}-\frac{1}{2}}\right \}.
$$
In what follows we consider the remaining cases where $\alpha$ and
$\beta$ are distinct and at least one of them is $0$ or $2$.

Consider the sequences $\bl =\{ d_n \}$, $\bl'= \{ d'_n \}$
and $\bl''= \{ d''_n \}$ defined by
$$d_n:=2^{-\frac{n}{\log n}}, \quad
\quad d'_n:=2^{-\nu n}, \quad \quad d''_n:=2^{-n \log n},$$
where $\nu>0$, and construct the Cantor sets $\Lambda=\Lambda
(\bl)$, $\Lambda'=\Lambda(\bl')$ and $\Lambda''=\Lambda(\bl'')$ as
in \S \ref{sec:bg}. By \lemref{Hdim},
$$\Hdim(\Lambda)=2, \quad \quad \Hdim(\Lambda')
=\frac{2}{\nu+1}, \quad \quad \Hdim(\Lambda'')=0.$$
We prove that the standard homeomorphisms between these three
Cantor sets {\it and} their inverses are all David maps;
this will prove the theorem. In view of Tukia's Theorem quoted in 
\S \ref{sec:intro}, it suffices to check that the sequence of
approximating homeomorphisms are David maps with uniform parameters
$(C,\alpha,K_0)$. In fact, the estimates below show that we can always 
take $C=\alpha=1$. \vs

$\bullet$ {\it Case 1. Mapping $\Lambda$ to $\Lambda'$.} Suppose
$\{ \varphi_n \}$ is the sequence of quasiconformal maps which
approximates the standard homeomorphism $\varphi$ from $\Lambda$ to
$\Lambda'$. To estimate the dilatation of $\varphi_n$, note that
\begin{equation}
\label{eqn:an}
a_n=2^{-(n+1)}(d_{n-1}-d_n) \asymp 2^{-n}
(2^{-\frac{n-1}{\log (n-1)}}-2^{-\frac{n}{\log n}}) \asymp
\frac{2^{-n-\frac{n}{\log n}}}{\log n}
\end{equation}
and
\begin{equation}
\label{eqn:a'n}
a'_n=2^{-(n+1)}(d'_{n-1}-d'_n) \asymp 2^{-n}
(2^{-\nu(n-1)}-2^{-\nu n}) \asymp 2^{-(\nu+1)n}.
\end{equation}
Hence
$$
\frac{a_n'\, d_{n}}{a_n\, d'_{n}} \asymp \frac{2^{-(\nu+1)n} \cdot
2^{-\frac{n}{\log n}}}{\frac{2^{-n-\frac{n}{\log n}}}{\log n} \cdot
2^{-\nu n}} \asymp \log n.
$$
It follows from (\ref{eqn:dil}) that there is a sequence
$1<K_1<K_2< \cdots < K_n < \cdots$ with $K_n \asymp \log n$
such that $\varphi_n$ is $K_n$-quasiconformal. Fix the index $n$
and a number $K>1$. Choose $j$ so that $K_j \leq K < K_{j+1}$. Then
\begin{align*}
\area \{ z: K_{\varphi_n}(z) > K \} \leq & \, \area \{ z:
K_{\varphi_n}(z) > K_j \} \vs \\
\leq & \, \area (\Lambda_j)=d_j^2=4^{-\frac{j}{\log j}}.
\end{align*}
Since $K \asymp K_j \asymp \log j$, we obtain 
$$\area \{ z: K_{\varphi_n}(z) > K \} \leq e^{-K},$$
provided that $K$ is bigger than some $K_0$ independent of $n$. 
It follows that the $\varphi_n$ are all $(1,1,K_0)$-David maps.

The inverse maps $\psi_n:=\varphi_n^{-1}$ are also
$K_n$-quasiconformal with the same dilatation $K_n \asymp \log n$
and they converge uniformly to $\psi:=\varphi^{-1}$. Moreover, if
$K_j \leq K < K_{j+1}$, then
\begin{align*}
\area \{ z: K_{\psi_n}(z) > K \} \leq & \, \area \{ z:
K_{\psi_n}(z) > K_j \} \vs \\
\leq & \, \area (\Lambda'_{j})=(d'_{j})^2=4^{-\nu \, j} \vs \\
\leq & \, e^{-K},
\end{align*}
provided that $K$ is bigger than some $K_0$ independent of $n$. 
It follows that the $\psi_n$ are all $(1,1,K_0)$-David maps. \vs

$\bullet$ {\it Case 2. Mapping $\Lambda'$ to $\Lambda''$.} The
argument here is quite similar to the previous case. We have
\begin{equation}
\label{eqn:a''n}
a''_n=2^{-(n+1)}(d''_{n-1}-d''_n) \asymp 2^{-n}
(2^{-(n-1)\log(n-1)}-2^{-n\log n}) \asymp 2^{-n-n\log n+\log n}
\end{equation}
Hence, using \eqref{eqn:a'n} and \eqref{eqn:a''n}, we obtain
$$
\frac{a''_n\, d'_{n}}{a'_n\, d''_{n}} \asymp
\frac{2^{-n-n\log n+\log n} \cdot 2^{-\nu n}}{2^{-(\nu+1)n} \cdot
2^{-n\log n}} \asymp 2^{\log n}.
$$
Let $\{ \varphi_n \}$ be the sequence of quasiconformal maps which
approximates the standard homeomorphism $\varphi$ from $\Lambda'$
to $\Lambda''$. It follows from (\ref{eqn:dil}) that there is
a sequence $1<K_1<K_2< \cdots < K_n < \cdots$ with $K_n \asymp
2^{\log n}$ such that $\varphi_n$ is $K_n$-quasiconformal. Fix the
index $n$ and a number $K>1$, and choose $j$ so that $K_j \leq K <
K_{j+1}$. Then $K \asymp K_j \asymp 2^{\log j}$ and
\begin{align*}
\area \{ z: K_{\varphi_n}(z) > K \} \leq & \, \area \{ z:
K_{\varphi_n}(z) > K_j \} \vs \\
\leq & \, \area (\Lambda'_{j}) = (d'_{j})^2=4^{-\nu j} \vs \\
\leq & \, e^{-K},
\end{align*}
provided that $K$ is bigger than some $K_0$ independent of $n$. 

The inverse maps $\psi_n:=\varphi_n^{-1}$ are $K_n$-quasiconformal
with $K_n \asymp 2^{\log n}$ and they converge uniformly to
$\psi:=\varphi^{-1}$. Moreover, if $K_j \leq K < K_{j+1}$, then
\begin{align*}
\area \{ z: K_{\psi_n}(z) > K \} \leq & \, \area \{ z:
K_{\psi_n}(z) > K_j \} \vs \\
\leq & \, \area (\Lambda''_{j})=(d''_{j})^2=4^{-j \log j} \vs \\
\leq & \, e^{-K},
\end{align*}
provided that $K$ is bigger than some $K_0$ independent of $n$. \vs

$\bullet$ {\it Case 3. Mapping $\Lambda$ to $\Lambda''$.} Using
(\ref{eqn:an}) and (\ref{eqn:a''n}), we obtain
$$
\frac{a''_n\, d_{n}}{a_n\, d''_{n}} \asymp
\frac{2^{-n-n\log n+\log n}\cdot 2^{-\frac{n}{\log n}}}
{\frac{2^{-n-\frac{n}{\log n}}}{\log n} \cdot 2^{-n \log n}} \asymp
2^{\log n}\ \log n= n^{\log 2} \ \log n.
$$
Let $\{ \varphi_n \}$ be the sequence of quasiconformal maps which
approximates the standard homeomorphism $\varphi$ from $\Lambda$ to
$\Lambda''$. It follows then from (\ref{eqn:dil}) that there is a
sequence $1<K_1<K_2< \cdots < K_n < \cdots$ with $K_n \asymp
n^{\log 2}\ \log n$ such that $\varphi_n$ is $K_n$-quasiconformal.
Fix $n$, let $K$ be sufficiently large, and choose $j$ so that $K_j
\leq K < K_{j+1}$. Then
\begin{align*}
\area \{ z: K_{\varphi_n}(z) > K \} \leq & \, \area \{ z:
K_{\varphi_n}(z) > K_j \} \vs \\
\leq & \, \area (\Lambda_{j}) = (d_{j})^2=4^{-\frac{j}{\log j}}.
\end{align*}
But $K \asymp K_j \asymp j^{\log 2}\ \log j$, so
$$\area \{ z: K_{\varphi_n}(z) > K \} \leq e^{-K},$$
provided that $K$ is bigger than some $K_0$ independent of $n$. 

The inverse maps $\psi_n:=\varphi_n^{-1}$ are $K_n$-quasiconformal
with $K_n \asymp n^{\log 2}\ \log n$ and they converge uniformly to
$\psi:=\varphi^{-1}$. Moreover, if $K_j \leq K < K_{j+1}$, then
\begin{align*}
\area \{ z: K_{\psi_n}(z) > K \} \leq & \, \area \{ z:
K_{\psi_n}(z) > K_j \} \vs \\
\leq & \, \area (\Lambda''_{j})=(d''_{j})^2=4^{-j \log j} \vs \\
\leq & \, e^{-K}, 
\end{align*}
provided that $K$ is bigger than some $K_0$ independent of $n$. \hfill
$\Box$

\section{Proof of Theorem B}
\label{sec:proofB}

The idea of the proof of Theorem B is to construct a David
map $\varphi: \CC \to \CC$ which sends a {\it linear}
Cantor set $\Sigma \subset [-\frac{1}{2},\frac{1}{2}]$ to a Cantor
set of the form $\Lambda(\bl)$ with dimension $2$. The image
$\varphi([-\frac{1}{2},\frac{1}{2}])$ will then be an
embedded arc of dimension $2$. Since the construction allows
$\varphi=\text{id}$ outside the square $[-\frac{1}{2},\frac{1}{2}]
\times [-\frac{1}{2},\frac{1}{2}]$, we can easily complete this arc
to a David circle.

\subsection*{A linear Cantor set}
Consider the closed unit square
$\Sigma_0:=[-\frac{1}{2},\frac{1}{2}] \times
[-\frac{1}{2},\frac{1}{2}]$ in the plane. We construct a nested
sequence $\{ \Sigma_n \}_{n \geq 0}$ of compact sets whose
intersection is a linear Cantor set. For $1\leq j \leq 4$, let
$f_j: \CC \to \CC$ be the affine contraction defined by
$$
f_j(z)=\frac{1}{8}z+\frac{2j-5}{8},
$$
and set
$$
\Sigma_n:= \bigcup_{j_1, \ldots , j_n} f_{j_1} \circ \cdots \circ
f_{j_n} (\Sigma_0),
$$
where the union is taken over all unordered $n$-tuples
$j_1,\ldots,j_n$ chosen from $\{ 1,2,3,4 \}$. It is easy to see
that $\Sigma_n$ is the disjoint union of $4^n$ closed squares of
side-length $8^{-n}$ with centers on $[-\frac{1}{2},\frac{1}{2}]$
and sides parallel to the coordinate axes (compare \figref{ds} left).
We define the Cantor set
$\Sigma$ as the intersection $\bigcap_{n=0}^{\infty} \Sigma_n$.
Evidently, $\Sigma$ is a subset of $[-\frac{1}{2},\frac{1}{2}]$
which has linear measure zero and Hausdorff dimension $2/3$.

\realfig{curve}{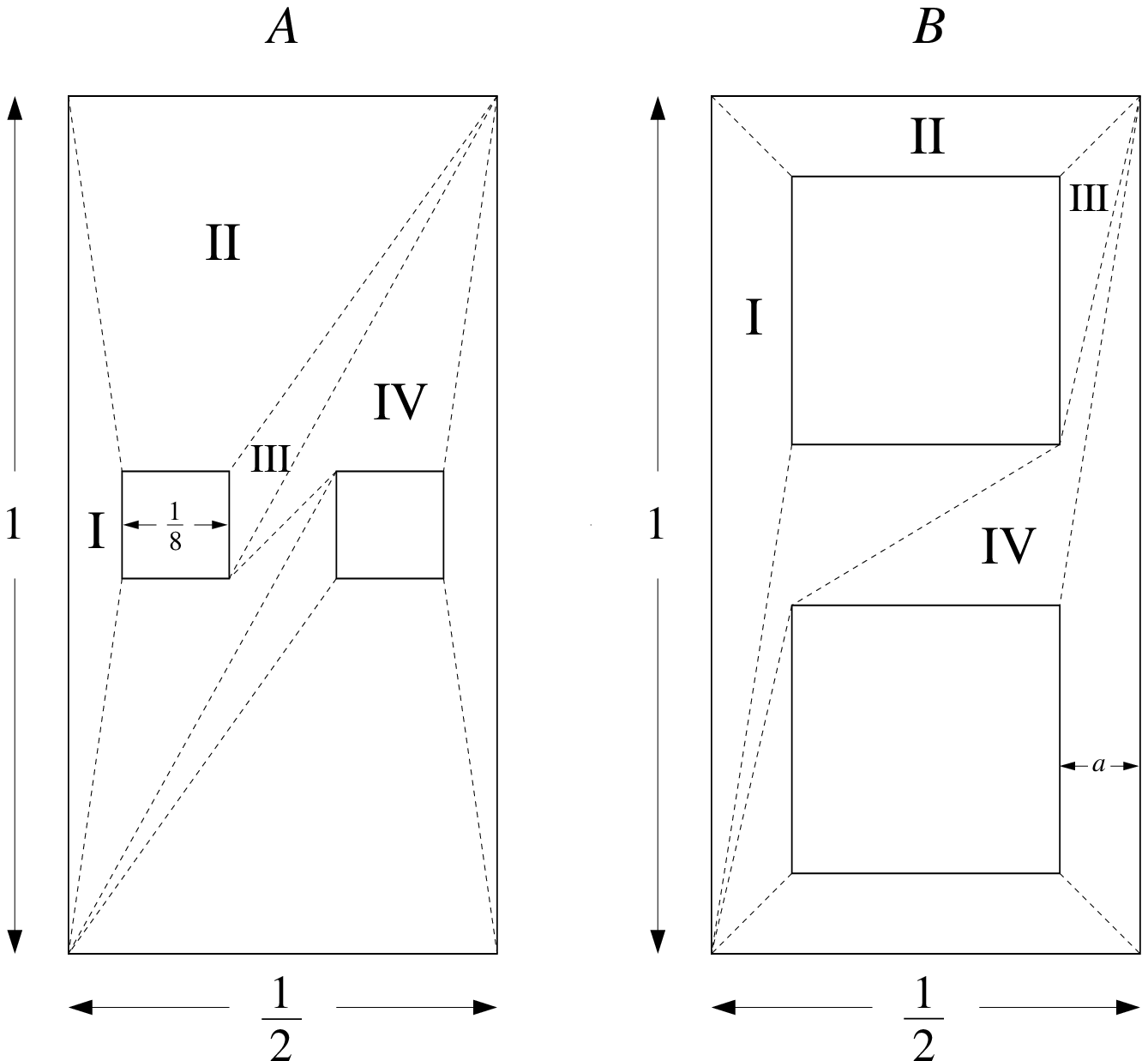}{{\sl Cell decompositions of $A$ and
$B$.}}{8cm}

\subsection*{A quasiconformal twist}
The proof of Theorem B depends on the following lemma which is a
triply-connected version of \lemref{annulus}. For simplicity we
denote by $S(p,r)$ the open square centered at $p$ whose
side-length is $r$.

\begin{lemma}
\label{twist}
Fix $0<a<1/5$ and let $A$ and $B$ be the closed triply-connected
sets defined by
\begin{align*}
A:= & \left( \left[ 0 , \frac{1}{2} \right] \times \left[
-\frac{1}{2}, \frac{1}{2} \right] \right) \sm \left( S
\left(\frac{1}{8},\frac{1}{8} \right)
\cup S \left( \frac{3}{8},\frac{1}{8} \right) \right) \vs \\
B:= & \left( \left[ 0 , \frac{1}{2} \right] \times \left[
-\frac{1}{2}, \frac{1}{2} \right] \right) \sm \left( S
\left(\frac{1+i}{4}, \frac{1}{2}-2a \right) \cup S \left(
\frac{1-i}{4},\frac{1}{2}-2a \right) \right)
\end{align*}
(see \figref{curve}). Let $\varphi: \bd A \to \bd B$ be a
homeomorphism which is the identity on the outer boundary component
and acts affinely on the inner boundary components, mapping $\bd
S(\frac{1}{8},\frac{1}{8})$ to $\bd S(\frac{1+i}{4},
\frac{1}{2}-2a)$ and $\bd S(\frac{3}{8},\frac{1}{8})$ to $\bd
S(\frac{1-i}{4}, \frac{1}{2}-2a)$, respecting the horizontal and
vertical sides. Then $\varphi$ can be extended to a
$K$-quasiconformal map $\varphi: A \to B$, with
$$
K \asymp \frac{1}{a} \, .
$$
\end{lemma}

\realfig{norm}{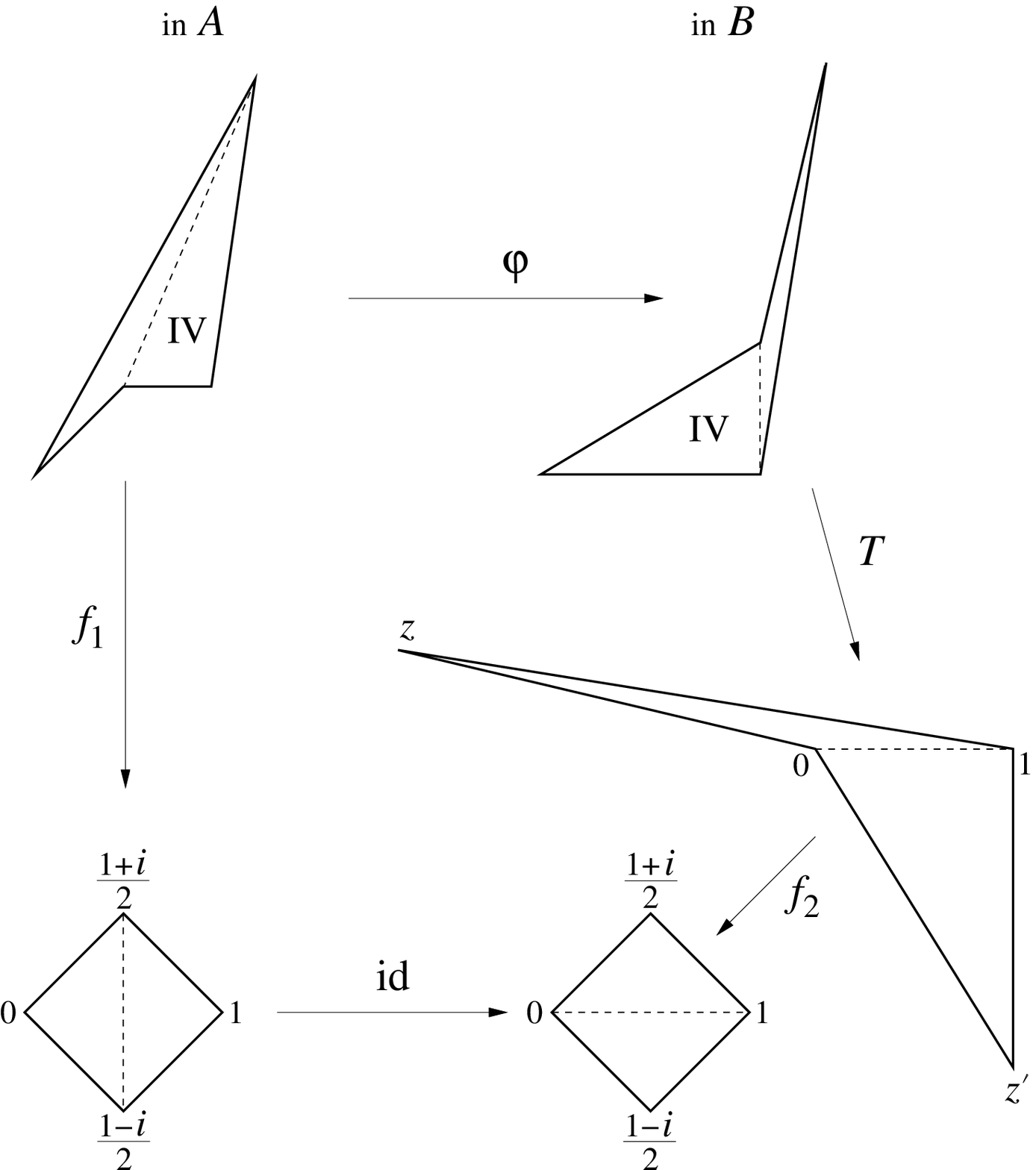}{{\sl Extending $\varphi$ between cells of
type IV.}}{8cm}

\begin{proof}
We consider the affine cell decompositions of $A$ and $B$ shown in
\figref{curve} and require $\varphi$ to map each cell in $A$ to its
corresponding cell in $B$ in a piecewise affine fashion. By
symmetry, it suffices to define $\varphi$ piecewise affinely
between the cells labeled I, II, III, and IV. We let $\varphi$ be
affine between the triangular cells III. On the cells I and II we
subdivide the trapezoids into two triangular cells and define
$\varphi$ to be affine on each of them. An easy computation based
on \eqref{eqn:KL} then shows that the dilatation of $\varphi$ on I,
II, and III is comparable to $1/a$.

It remains to define $\varphi$ between the cells IV and estimate
its dilatation. Note that the cell IV in $A$ has bounded geometry,
so there is a $K_1 \asymp 1$ and a piecewise affine $K_1$-quasiconformal
map $f_1$ from this cell to the square with vertices
$0, 1, (1+i)/2, (1-i)/2$ which maps
the horizontal edge of this cell to the segment from $(1-i)/2$ to $1$
(see \figref{norm}). The cell IV in $B$, after a conformal
change of coordinates $T$, becomes the $4$-gon with vertices
$$0,\ 1,\ z:=-\frac{1-2a}{4a}+\frac{i}{2},\ z':=1-\frac{(1-4a)i}{4a}.$$
Let $f_2$ be the piecewise affine map on this $4$-gon which maps
the triangle $\Delta(0,1,z)$ to $\Delta(0,1,(1+i)/2)$ and
the triangle $\Delta(0,1,z')$ to $\Delta(0,1,(1-i)/2)$
(see \figref{norm}). Then a brief calculation based on \eqref{eqn:KL}
shows that $f_2$ is $K_2$-quasiconformal, with $K_2 \asymp 1/a$.
The map $\varphi$ can then be defined by
$T^{-1} \circ f_2^{-1} \circ f_1$, whose
dilatation $K_1K_2$ is clearly comparable to $1/a$.
\end{proof}

We are now ready to prove Theorem B cited in \S \ref{sec:intro}. \vs \\
{\it Proof of Theorem B.} Consider the Cantor set
$\Sigma=\bigcap_{n=0}^{\infty} \Sigma_n$ constructed above and the
Cantor set $\Lambda=\Lambda(\bl)=\bigcap_{n=0}^{\infty} \Lambda_n$
constructed in \S \ref{sec:bg}, where
$\bl = \{ d_n \}$ is defined by $d_n:=2^{-\sqrt{n}}$. It follows from
\lemref{Hdim} that $\Hdim(\Lambda)=2$.

We construct a David map $\varphi: \CC \to \CC$, identity
outside the square $[-\frac{1}{2}, \frac{1}{2}] \times
[-\frac{1}{2}, \frac{1}{2}]$, with the property
$\varphi(\Sigma)=\Lambda$. Then the embedded arc
$\varphi([-\frac{1}{2}, \frac{1}{2}])$ contains $\Lambda$ and hence
has dimension $2$. By pre-composing $\varphi$ with an appropriate
quasiconformal map, we obtain a David map sending the
round circle to a Jordan curve of dimension $2$.

The map $\varphi$ will be the uniform limit of a sequence of quasiconformal
maps $\varphi_n: \CC \to \CC$ with $\varphi_n(\Sigma_n)=\Lambda_n$,
defined inductively as follows. Let $\varphi_0$ be the identity map on
$\CC$. To define $\varphi_1$, set $\varphi_1=\varphi_0$ on $\CC \sm
\Sigma_0$ and map each of the four squares in $\Sigma_1$ affinely to the
``corresponding'' square in $\Lambda_1$. Here ``corresponding'' means that
the squares in $\Sigma_0$, from left to right, map respectively
to the north west, south west, north east and south east squares
in $\Lambda_1$ (compare \figref{ds}). The remaining set
$\Sigma_0 \sm \Sigma_1$ is the union of two triply-connected regions,
on the boundary of which $\varphi_1$ can be defined affinely, so we can extend
$\varphi_1$ to each such region as in \lemref{twist}.

In general, suppose $\varphi_{n-1}$ is constructed for some $n \geq 2$ and
that it maps each square in $\Sigma_{n-1}$ affinely to a square
in $\Lambda_{n-1}$. Define $\varphi_n=\varphi_{n-1}$ on $\CC \sm
\Sigma_{n-1}$ and let $\varphi_n$ map each square in $\Sigma_n$
affinely to the ``corresponding'' square in $\Lambda_n$ in the above sense.
The remaining set $\Sigma_{n-1} \sm \Sigma_n$ is the union of $2^{2n-1}$
triply-connected regions on the boundary of which $\varphi_n$ can be
defined affinely.
By rescaling each such region in $\Sigma_{n-1} \sm \Sigma_n$ by a factor
$8^{n-1}$ and the corresponding region in $\Lambda_{n-1} \sm \Lambda_n$ by a
factor $2^{n-1}/d_{n-1}$, we are in the situation of \lemref{twist},
so we can extend $\varphi_n$ in a piecewise affine fashion as in that lemma,
and the dilatation of the resulting extension will be comparable to
$$\begin{array}{rl}
\ds{\frac{d_{n-1}}{2^{n-1}a_n}} & \ds{= \frac{d_{n-1}}{2^{n-1}\cdot
2^{-(n+1)}(d_{n-1}-d_n)}} \vs \\
& \ds{= \frac{2^{-\sqrt{n-1}}}{2^{n-1}\cdot 2^{-(n+1)}
(2^{-\sqrt{n-1}}-2^{-\sqrt{n}})}} \vs \\
& \asymp \sqrt{n}.
\end{array}$$

\realfig{ds}{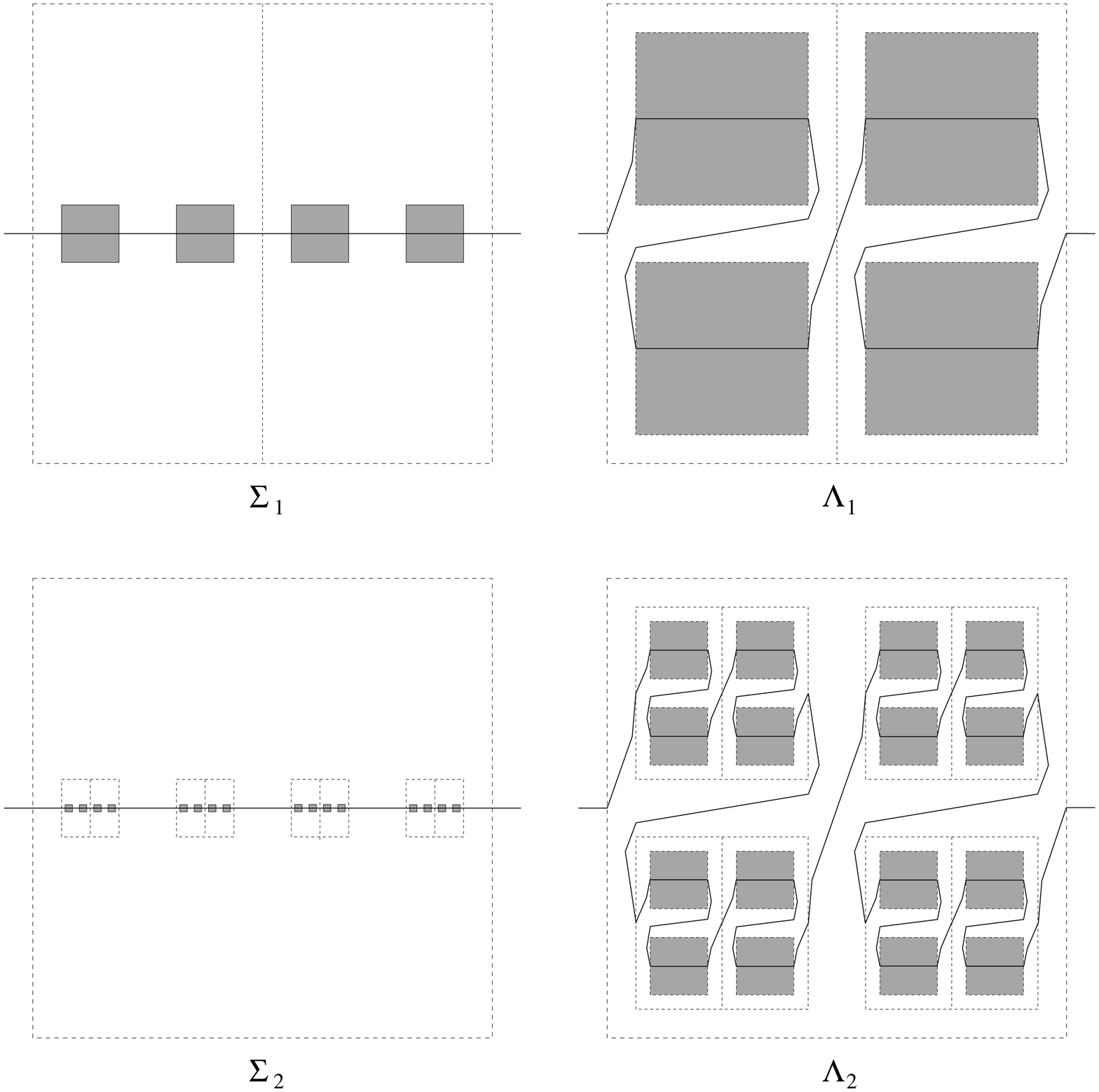}{{\sl First two steps in the construction
of the map $\varphi$. The solid arcs on the right are
$\varphi_n(\RR)$ for $n=1,2$.}}{12cm}

The sequence $\{ \varphi_n \}$ obtained this way has the following properties:
\begin{enumerate}
\item[(i)]
$\varphi_n=\varphi_{n-1}$ on $\CC \sm \Sigma_{n-1}$. \vs
\item[(ii)]
$\varphi_n$ maps each square in $\Sigma_n$ affinely to the
corresponding square in $\Lambda_n$. \vs
\item[(iii)]
$\varphi_n$ is $K_n$-quasiconformal, with $K_n \asymp \sqrt{n}$. \vs
\end{enumerate}
Evidently, $\varphi:=\lim_{n \to \infty} \varphi_n$ is a
homeomorphism which agrees with $\varphi_n$ on $\CC \sm \Sigma_n$
for every $n$ and satisfies $\varphi (\Sigma)=\Lambda$.

To check that $\varphi$ is a David map, choose a sequence
$1 < K_1 < K_2 < \cdots < K_n < \cdots$ with $K_n \asymp \sqrt{n}$
such that $\varphi_n$ is $K_n$-quasiconformal. Fix some $n$, let $K>1$, 
and choose $j$ such that $K_j \leq K < K_{j+1}$. Then
\begin{align*}
\area \{ z: K_{\varphi_n}(z) > K \} \leq & \, \area \{ z:
K_{\varphi_n}(z) > K_j \} \vs \\
\leq & \, \area (\Sigma_j)=2^{-4j}.
\end{align*}
Since $K \asymp K_j \asymp \sqrt{j}$, we have
$$\area \{ z: K_{\varphi_n}(z) > K \} \leq e^{-K},$$
provided that $K$ is bigger than some $K_0$ independent of $n$. 
It follows that the $\varphi_n$ are all $(1,1,K_0)$-David maps. 
By Tukia's Theorem in \S \ref{sec:intro},
we conclude that $\varphi=\lim_{n \to \infty}
\varphi_n$ is a David map. \hfill $\Box$

\subsection*{Removability of David circles}

A compact set $\Gamma \subset \CC$ is called {\it
(quasi)conformally removable} if every homeomorphism $\varphi: \CC
\to \CC$ which is (quasi)conformal off $\Gamma$ is (quasi)conformal
in $\CC$. It is well-known that conformal and quasiconformal
removability are identical notions.

Every set of $\sigma$-finite $1$-dimensional Hausdorff measure,
such as a rectifiable curve, is removable. Quasiarcs and
quasicircles provide examples of removable sets which can have any
dimension in the interval $[1,2)$. One can even construct removable
sets of dimension $2$: the Cartesian product of two linear Cantor
sets with zero length and dimension $1$ is such a set.

At the other extreme, sets of positive area are never removable,
as can be seen by an easy application of the measurable Riemann
mapping theorem. Also, there exist non-removable sets of Hausdorff
dimension $1$ (see for example \cite{K}).

To add an item to the above list of examples, we show that David circles are
removable, which, combined with Theorem B, proves that {\it there exist
removable Jordan curves of Hausdorff dimension $2$}. First we need the
following simple lemma on David maps (compare \cite{PZ})
whose analogue in the quasiconformal case is standard.

\begin{lemma}
\label{sewing}
Suppose $\varphi: \CC \to \CC$ is a homeomorphism whose restrictions to
$\DD$ and $\CC \sm \ov{\DD}$ are David. Then $\varphi$ itself is a David
map.
\end{lemma}

\begin{proof}
The complex dilatation $\mu=\mu_{\varphi}$ is defined almost
everywhere in $\CC$ and satisfies an exponential condition of the
form (\ref{eqn:David1}) in $\DD$ and in $\CC \sm \ov{\DD}$ 
(by making $C$ bigger and $t$ and $\ve_0$ smaller if necessary, we
can assume that the same constants $(C,t,\ve_0)$ work for both $\DD$
and $\CC \sm \ov{\DD}$). So to prove the lemma, we need only show that $\varphi
\in W^{1,1}_{\text{loc}}(\CC)$.

On every compact subset of $\CC \sm \Sen$, the ordinary partial derivatives
$\bd \varphi$ and $\ov{\bd}\varphi$ exist almost everywhere, are
integrable, and coincide with the distributional partial
derivatives of $\varphi$. We check that $\bd \varphi$, and hence
$\ov{\bd}\varphi=\mu \cdot \bd \varphi$, is locally integrable near
the unit circle $\Sen$.

Let $D$ be any small disk centered on $\Sen$. We have
$$
|\bd \varphi|^2 = \frac{J_{\varphi}}{1-|\mu|^2} \leq
\frac{J_{\varphi}}{1-|\mu|},
$$
so that
\begin{equation}
\label{eqn:ho}
|\bd \varphi| \leq (J_{\varphi})^{\frac{1}{2}} \cdot (1-|\mu|)^
{-\frac{1}{2}}.
\end{equation}
Now $J_{\varphi} \in L^1(D)$ since $\int_D J_{\varphi} \leq
\area(\varphi(D)) < +\infty$, and $(1-|\mu|)^{-1} \in L^1(D)$
because of the exponential condition (\ref{eqn:David1}). It follows from 
H\"{o}lder inequality applied to \eqref{eqn:ho} that 
$\bd \varphi \in L^1(D)$.
\end{proof}

\begin{theorem}
David circles are (quasi)conformally removable.
\end{theorem}

\begin{proof}
Let $\varphi: \CC \to \CC$ be a David map and
$\Gamma=\varphi(\Sen)$. Let $f: \CC \to \CC$ be a homeomorphism
which is conformal in $\CC \sm \Gamma$. Then the homeomorphism
$f \circ \varphi$ is David in $\DD$ and in $\CC \sm \ov{\DD}$.
By \lemref{sewing}, $f \circ \varphi : \CC \to \CC$ is a David
map. Since $\mu_{f \circ \varphi}=\mu_{\varphi}$ almost
everywhere, it follows from the uniqueness part of David's theorem
\cite{David} that $f$ must be conformal in $\CC$.
\end{proof}


\begin{thebibliography}{*****}
\bibitem [\bf{A}]{Ahlfors} L.~Ahlfors, {\it Lectures on quasiconformal
mappings}, Van Nostrand, 1966.
\bibitem [\bf{AB}]{AB} L.~Ahlfors and L.~Bers,
{\it Riemann mapping's theorem for variable metrics},
Annals of Math. {\bf 72} (1960) 385-404.
\bibitem [\bf{As}]{Astala} K.~Astala, {\it Area distortion of
quasiconformal mappings}, Acta Math. {\bf 173} (1994) 37-60.
\bibitem [\bf{D}]{David} G.~David,
{\it Solutions de l'equation de Beltrami avec $\| \mu \|_{\infty} = 1$},
Ann. Acad. Sci. Fenn. Ser. A I Math. {\bf 13} (1988) 25-70.
\bibitem [\bf{GV}]{GV} F.~Gehring, J.~V\"ais\"al\"a, {\it Hausdorff
dimension and quasiconformal mappings}, J. London Math. Soc. {\bf 6}
(1973) 504-512.
\bibitem [\bf{GJ}]{GJ} J.~Graczyk and P.~Jones, {\it Dimension of
the boundary of quasiconformal Siegel disks}, Invent. Math. {\bf 148}
(2002), 465-493.
\bibitem [\bf{K}]{K} R.~Kaufman, {\it Fourier-Stieltjes coefficients 
and continuation of functions}, Ann. Acad. Sci. Fenn. Ser. A I
Math. {\bf 9} (1984) 27-31.
\bibitem [\bf{LV}]{LV} O.~Lehto, K.~Virtanen, {\it Quasiconformal
mappings in the plane}, Springer-Verlag, 1973.
\bibitem [\bf{M}]{Mattila} P.~Mattila, {\it Geometry of Sets and
Measures in Euclidean spaces}, Cambridge University Press, 1995.
\bibitem [\bf{P}]{Petersen} C.~L.~Petersen, {\it Local connectivity of
some Julia sets containing a circle with an irrational rotation}, Acta
Math. {\bf 177} (1996) 163-224.
\bibitem [\bf{PZ}]{PZ} C.~L.~Petersen, S.~ Zakeri, {\it On the Julia
set of a typical quadratic polynomial with a Siegel disk}, Stony
Brook IMS preprint 2000/6, to appear in Annals of Math.
\bibitem [\bf{T}]{Tukia} P.~Tukia,
{\it Compactness properties of $\mu$-homeomorphisms},
Ann. Acad. Sci. Fenn. Ser. A I Math. {\bf 16} (1991) 47-69.

\end{thebibliography}
\end{document}